\documentclass[reqno]{amsart}

\usepackage{amsmath,amssymb,amsthm,amsfonts,mathrsfs}
\usepackage{fullpage}
\usepackage{xcolor}
\usepackage[colorlinks=true,
linkcolor=blue,
anchorcolor=blue,
citecolor=red
]{hyperref}
\allowdisplaybreaks 
\newtheorem{theorem}{Theorem}[section]
\newtheorem{lemma}[theorem]{Lemma}

\newtheorem{proposition}[theorem]{Proposition}

\newtheorem*{conjecture*}{Conjecture}
\theoremstyle{definition}
\newtheorem{definition}[theorem]{Definition}

\theoremstyle{remark}

\newtheorem*{remark*}{remark}

\usepackage{colonequals}

\author{Runbo Li}
\address{The High School Affiliated to Renmin University of China International Curriculum Center, Beijing, China}
\email{runbo.li.carey@gmail.com}

\makeatletter
\@namedef{subjclassname@2020}{\textup{}2020 Mathematics Subject Classification}
\makeatother

\title[]{On the largest prime factor of quadratic polynomials}

\begin{document}
	
\begin{abstract}
Let $x$ denote a sufficiently large integer. We show that the recent result of Grimmelt and Merikoski actually yields the largest prime factor of $n^2 +1$ is greater than $x^{1.317}$ infinitely often. As an application, we give a new upper bound for the number of integers $n \leqslant x$ which $n^2 +1$ has a primitive divisor.
\end{abstract}

\maketitle

\tableofcontents

\section{Introduction}
Let $x, n$ denote sufficiently large integers, $p$ denote a prime number, $P_{r}$ denote an integer with at most $r$ prime factors counted with multiplicity, and let $f$ be an irreducible polynomial with degree $g$. It's conjectured that there are infinitely many $n$ such that $f(n)$ is prime. The simplest case is $g=1$, which is the famous Dirichlet's theorem proved more than 100 years ago. However, for $g \geqslant 2$, this conjecture is still open.

For the second simplest case $g=2$, there are several ways to attack this conjecture. One way is to relax the number of prime factors of $f(n)$, and the best result in this way is due to Iwaniec \cite{IwaniecP2}. Building on the previous work of Richert \cite{RichertP3}, he showed that for any irreducible polynomial $f(n)=an^2 +bn+c$ with $a>0$ and $c \equiv 1 (\bmod 2)$, there are infinitely many $x$ such that $f(x)$ is a $P_2$.

Another possible way is to consider the largest prime factor of $f(n)$. Let $P^{+}(x)$ denote the largest prime factor of $x$, then we hope to show that the largest prime factor of $f(n)$ is greater than $n^g$ for infinitely many integers $n$. For general polynomials, the best result is due to Tenenbaum \cite{Tenenbaum}, where he showed that for some $0<t<2-\log 4$, the largest prime factor of $f(n)$ is greater than $n \exp( (\log n)^{t} )$ for infinitely many integers $n$. However, it's rather difficult to prove the same thing holds for $n^{1+\varepsilon}$ even for a small $\varepsilon$.

For the special case $f(n)=n^2 +1$, the progress is far more than the general case. In 1967, Hooley \cite{Hooley} first proved the largest prime factor of $n^2 +1$ is greater than $n^{1.10014}$ for infinitely many integers $n$ by using the Weil bound for Kloosterman sums. By applying their new bounds for multilinear forms of Kloosterman sums, Deshouillers and Iwaniec \cite{DI82} showed in 1982 that the largest prime factor of $n^2 +1$ is greater than $n^{1.202468}$ infinitely often. In 2020, de la Bretèche and Drappeau \cite{BD2019} improved the exponent to $1.2182$ by making use of the result of Kim and Sarnak \cite{KimSarnak}. In 2023, Merikoski \cite{Merikoski} proved a new bilinear estimate and used Harman's sieve to get the exponent $1.279$. This is the first attempt of using Harman's sieve on this problem. In 2024, Pascadi \cite{Pascadi} optimized the exponent to $1.3$ by inserting his new arithmetic information. Recently, using a different approach to obtain arithmetical information, Grimmelt and Merikoski \cite{GrimmeltMerikoski} got $1.312$, which is the value that previously obtained by Merikoski \cite{Merikoski} under the Selberg eigenvalue conjecture. In the present paper, we shall use the exactly same sieve argument as in \cite{GrimmeltMerikoski} and illustrate that this exponent can be further improved to $1.317$.
\begin{theorem}\label{t1}
Let $(\frac{x}{p})$ denotes the Legendre symbol. There exists some small $\varepsilon > 0$ such that the following holds for all $X > 1/\varepsilon$. Let $1 \leqslant h \leqslant X^{1+\varepsilon}$ be square--free and $1 \leqslant a \leqslant X^{\varepsilon}$ with $(a,h)=1$. Suppose that
$$
\left| \sum_{p \leqslant Y}\frac{\log p}{p} \left(\frac{-ah}{p} \right) \right| \leqslant \varepsilon \log Y
$$
for any $X^{\varepsilon} < Y \leqslant X^2$. Then there exists $n \in [X, 2X]$ such that the largest prime factor of $a n^2 +h$ is greater than $n^{1.317}$. Specially, the largest prime factor of $n^2 + 1$ is greater than $n^{1.317}$ infinitely often.
\end{theorem}

As an application of our Theorem~\ref{t1}, we consider the polynomial $n^2 +1$ with a primitive divisor.
\begin{definition}
Let $\left(A_{n}\right)$ denote a sequence with integer terms. We say an integer $d>1$ is a primitive divisor of $A_{n}$ if $d \mid A_{n}$ and $\left(d, A_{m}\right)=1$ for all non-zero terms $A_{m}$ with $m<n$.
\end{definition}
\begin{proposition}
For all $n > 1$, the term $n^2 +1$ has a primitive divisor if and only if $P^{+}(n^2 +1) > 2n$. For all $n > 1$, if $n^2 +1$ has a primitive divisor then that primitive divisor is a prime and it is unique.
\end{proposition}
\noindent Contrary to the previous works on the lower bounds for the largest prime factor, a result due to Schinzel \cite{Schinzel} showed that for any $\varepsilon >0$, the largest prime factor of $n^2 +1$ is less than $n^{\varepsilon}$ infinitely often. In fact, from his result we can easily get the following.
\begin{theorem}([\cite{EverestHarman}, Theorem 1.2]).
The polynomial $n^2 +1$ does not have a primitive divisor for infinitely many terms.
\end{theorem}
\noindent We are interested in finding good upper and lower bound for the number of terms $n^2 +1$ with a primitive divisor. We define
$$
\rho(x)= \left| \left\{ n \leqslant x: n^{2}+1 \text { has a primitive divisor}\right\} \right|.
$$
Then we have the following simple upper bound
\begin{equation}
\rho(x) < x- \frac{C x}{\log x} 
\end{equation}
for some constant $C>0$. In \cite{ESTW} the following stronger result is mentioned.
\begin{equation}
\rho(x) < x- \frac{x \log \log x}{\log x} .
\end{equation}
In \cite{EverestHarman}, Everest and Harman first proved a lower bound with positive density and a better upper bound for $\rho(x)$. More precisely, they got the following bounds:
\begin{theorem}([\cite{EverestHarman}, Theorem 1.4]).
We have
$$
0.5324 x < \rho(x) < 0.905 x.
$$
\end{theorem}
\noindent They also conjectured the asymptotic $\rho(x) \sim (\log 2) x$ in their paper. In 2024, Harman \cite{Harman2024} used Merikoski's work on the largest prime factor of $n^2 +1$ and sharpened the upper and lower bounds for $\rho(x)$.
\begin{theorem}([\cite{Harman2024}, Theorem 5.5]).
We have
$$
0.5377 x < \rho(x) < 0.86 x.
$$
\end{theorem}
\noindent In the same year, Li \cite{LRBPrimitive} further improved the upper bound for $\rho(x)$ by using Pascadi's work.
\begin{theorem}([\cite{LRBPrimitive}, Theorem 1.6]).
We have
$$
\rho(x) < 0.847 x.
$$
\end{theorem}
\noindent Mine \cite{Mine} got a better lower bound for $\rho(x)$.
\begin{theorem}([\cite{Mine}, Theorem 1.3]).
We have
$$
\rho(x) > 0.543 x.
$$
\end{theorem}

In the present paper, we use the same sieve argument as in \cite{LRBPrimitive} and a recent result of Grimmelt and Merikoski to give a better upper bound for $\rho(x)$.
\begin{theorem}\label{t2}
We have
$$
\rho(x) < 0.838 x.
$$
\end{theorem}

\section{Merikoski's sieve decompositions}
Let $\varepsilon$ denote a sufficient small positive number and $P_{x}$ denote the largest prime factor of $\prod_{x \leqslant n \leqslant 2 x}\left(n^{2}+1\right)$. In this section we briefly introduce Grimmelt and Merikoski's work on finding a lower bound for $P_{x}$. Let $b(x)$ denote a non-nagative $C^{\infty}$--smooth function supported on $[x, 2x]$ and its derivatives satisfy $b^{(j)}(x) \ll x^{-j}$ for all $j \geqslant 0$. We define
$$
\left|\mathcal{A}_{d}\right|:=\sum_{n^{2}+1 \equiv 0 (\bmod d)} b(n) \quad \text{and} \quad X:=\int b(x) d x.
$$
Then by the method of Chebyshev--Hooley and the discussion in \cite{Merikoski}, we only need to find an upper bound for 
\begin{equation}
S(x):= \sum_{x<p \leqslant P_{x}}\left|\mathcal{A}_{p}\right| \log p = X \log x+O(x)
\end{equation}
with a constant less than 1. By a smooth dyadic partition we have
\begin{equation}
S(x)= \sum_{\substack{x \leqslant P \leqslant P_{x} \\ P=2^{j} x}} S(x, P)+O(x),
\end{equation}
where
\begin{equation}
S(x, P)=\sum_{P \leqslant p \leqslant 4 P} \psi_{P}(p)\left|\mathcal{A}_{p}\right| \log p
\end{equation}
for some $C^{\infty}$--smooth functions $\psi_{P}$ supported on $[P, 4 P]$ satisfying $\psi_{P}^{(l)}(x) \ll  P^{-l}$ for all $l \geqslant 0$.

In \cite{GrimmeltMerikoski}, Grimmelt and Merikoski proved the following upper bound for $S(x)$ with $P_{x} = x^{1.312}$ by using Harman's sieve method together with their new arithmetic information.

\begin{lemma}\label{l21}(See \cite{Merikoski}).
We have
\begin{align}
\nonumber \sum_{\substack{x \leqslant P \leqslant x^{1.312} \\ P=2^{j} x}} S(x, P) \leqslant&\ \left( G_0 + G_1 + G_2 + G_3 + G_4 + G_5 - G_6 + G_7 \right) X \log x \\
\nonumber <&\ 0.998 X \log x,
\end{align}
where
\begin{align}
\nonumber G_{0} & = \int_{1}^{\frac{7}{6}} 1 d \alpha = \frac{1}{6}, \\
\nonumber G_{1} & = \int_{1}^{\frac{17}{16}} \int_{\sigma(\alpha)}^{\alpha-2 \sigma(\alpha)} \alpha \frac{\omega\left(\frac{\alpha-\beta}{\beta}\right)}{\beta^2} d \beta d \alpha + \int_{1}^{\frac{17}{16}} \int_{\xi(\alpha)}^{\frac{\alpha}{2}} \alpha \frac{\omega\left(\frac{\alpha-\beta}{\beta}\right)}{\beta^2} d \beta d \alpha < 0.0287, \\
\nonumber G_{2} & = \int_{\frac{17}{16}}^{\frac{8}{7}} \int_{\sigma(\alpha)}^{\frac{\alpha}{2}} \alpha \frac{\omega\left(\frac{\alpha-\beta}{\beta}\right)}{\beta^2} d \beta d \alpha < 0.08622, \\
\nonumber G_{3} & = \int_{\frac{8}{7}}^{\frac{7}{6}} \int_{\sigma(\alpha)}^{\frac{\alpha}{2}} \alpha \frac{\omega\left(\frac{\alpha-\beta}{\beta}\right)}{\beta^2} d \beta d \alpha < 0.03107, \\
\nonumber G_{4} & =\int_{\frac{8}{7}}^{\frac{7}{6}} \int_{\sigma(\alpha)-\alpha+1}^{\alpha-1} \int_{\sigma(\alpha)-\alpha+1}^{\beta_1} \int_{\sigma(\alpha)-\alpha+1}^{\beta_2}
f_4\left(\alpha, \beta_1, \beta_2, \beta_3\right) \alpha \frac{\omega\left(\frac{\alpha-\beta_1 -\beta_2 -\beta_3}{\beta_3}\right)}{\beta_1 \beta_2 \beta_3^2} d \beta_3 d \beta_2 d \beta_1 d \alpha < 0.00011, \\
\nonumber G_{5} & = 4 \int_{\frac{7}{6}}^{\frac{5}{4}} \alpha d \alpha = \frac{29}{72}, \\
\nonumber G_{6} & = \int_{\frac{7}{6}}^{\frac{5}{4}} \int_{\alpha-1}^{\sigma(\alpha)} \alpha \frac{\omega\left(\frac{\alpha-\beta}{\beta}\right)}{\beta^2} d \beta d \alpha > 0.035631, \\
\nonumber G_{7} & = 4 \int_{\frac{5}{4}}^{1.312} \alpha d \alpha < 0.31769,
\end{align}
where
\begin{equation}
\sigma(\alpha) := \frac{2 -\alpha}{3}, \qquad \xi(\alpha) = \frac{3}{2} -\alpha,
\end{equation}
$f_4$ denotes the characteristic function of the set 
$$
\left\{ \beta_1 + \beta_2, \beta_1 + \beta_3, \beta_2 + \beta_3, \beta_1 + \beta_2 + \beta_3 \notin [\alpha-1, \sigma(\alpha)] \right\},
$$
and $\omega(u)$ denotes the Buchstab function determined by the following differential-difference equation
\begin{align*}
\begin{cases}
\omega(u)=\frac{1}{u}, & \quad 1 \leqslant u \leqslant 2, \\
(u \omega(u))^{\prime}= \omega(u-1), & \quad u \geqslant 2 .
\end{cases}
\end{align*}
\end{lemma}

However, their bounds for those integrals are not very accurate. Using Mathematica 14, we can get the following better bounds. We remark that for $G_{6}$ the new lower bound gives a $67\%$ improvement over the bound mentioned in \cite{Merikoski}.

\begin{lemma}\label{l22}
For $G_i$ ($0 \leqslant i \leqslant 6$) defined in Lemma~\ref{l21}, we have
$$
G_0 = \frac{1}{6}, \quad G_1 < 0.028611 (0.0287), \quad G_2 < 0.086062 (0.08622), \quad G_3 < 0.030992 (0.03107),
$$
$$
G_4 < 0.0001 (00.00011), \quad G_5 = \frac{29}{72}, \quad G_6 > 0.059841 (0.035631).
$$
Moreover, with these new bounds for $G_i$ ($0 \leqslant i \leqslant 6$) we have
$$
G_0 + G_1 + G_2 + G_3 + G_4 + G_5 - G_6 + 4 \int_{\frac{5}{4}}^{1.317} \alpha d \alpha < 0.9993.
$$
\end{lemma}
By Lemma~\ref{l22} and the same arguments as in \cite{GrimmeltMerikoski}, we complete the proof of Theorem~\ref{t1}.

\section{Proof of Theorem 1.9}
Let $V(u)$ denote an infinitely differentiable non-negative function such that
\begin{align*}
V(u)
\begin{cases}
<2, \quad & 1 < u < 2,\\
=0, \quad & u \leqslant 1 \text{ or } u \geqslant 2,
\end{cases}
\end{align*} 
with
$$
\frac{d^r V(u)}{d u^r} \ll 1 \quad \text{and} \quad \int_{\mathbb{R}} V(u) d u=1.
$$
By the discussion in \cite{EverestHarman} and \cite{Harman2024}, we wish to get an upper bound for sum of $\sum_{p \mid k^{2}+1} V(k/x)$ of the form
$$
\sum_{1 \leqslant p x^{-\alpha} \leqslant e} \sum_{p \mid k^{2}+1} V\left(\frac{k}{x}\right) \leqslant K(\alpha)(1+o(1)) \frac{X}{\log x}
$$
where $K(\alpha)$ is the sum of sieve theoretical functions related to the sieve decomposition on the problem of the largest prime factor of $n^2 +1$. This requires us to prove that for some $\tau$, we have
$$
\int_{1}^{\tau} \alpha K(\alpha) d \alpha < 1.
$$

By Lemma~\ref{l22} we can take $\tau =1.317$, and $K(\alpha)$ is defined as the piecewise function in Section 2. Combining this with the bound proved in \cite{EverestHarman}, we have
\begin{align}
\nonumber \rho(x) \leqslant&\ (1+o(1))x \int_{1}^{1.317} K(\alpha) d \alpha \\
\leqslant&\ \left( G^{\prime}_0 + G^{\prime}_1 + G^{\prime}_2 + G^{\prime}_3 + G^{\prime}_4 + G^{\prime}_5 - G^{\prime}_6 + G^{\prime}_7 \right) x ,
\end{align}
where
\begin{align}
\nonumber G^{\prime}_{0} & = \int_{1}^{\frac{7}{6}} \frac{1}{\alpha} d \alpha = \log \frac{7}{6} < 0.154151, \\
\nonumber G^{\prime}_{1} & = \int_{1}^{\frac{17}{16}} \int_{\sigma(\alpha)}^{\alpha-2 \sigma(\alpha)}  \frac{\omega\left(\frac{\alpha-\beta}{\beta}\right)}{\beta^2} d \beta d \alpha + \int_{1}^{\frac{17}{16}} \int_{\xi(\alpha)}^{\frac{\alpha}{2}}  \frac{\omega\left(\frac{\alpha-\beta}{\beta}\right)}{\beta^2} d \beta d \alpha < 0.027475, \\
\nonumber G^{\prime}_{2} & = \int_{\frac{17}{16}}^{\frac{8}{7}} \int_{\sigma(\alpha)}^{\frac{\alpha}{2}}  \frac{\omega\left(\frac{\alpha-\beta}{\beta}\right)}{\beta^2} d \beta d \alpha < 0.077933, \\
\nonumber G^{\prime}_{3} & = \int_{\frac{8}{7}}^{\frac{7}{6}} \int_{\sigma(\alpha)}^{\frac{\alpha}{2}}  \frac{\omega\left(\frac{\alpha-\beta}{\beta}\right)}{\beta^2} d \beta d \alpha < 0.026835, \\
\nonumber G^{\prime}_{4} & =\int_{\frac{8}{7}}^{\frac{7}{6}} \int_{\sigma(\alpha)-\alpha+1}^{\alpha-1} \int_{\sigma(\alpha)-\alpha+1}^{\beta_1} \int_{\sigma(\alpha)-\alpha+1}^{\beta_2}
f_4\left(\alpha, \beta_1, \beta_2, \beta_3\right)  \frac{\omega\left(\frac{\alpha-\beta_1 -\beta_2 -\beta_3}{\beta_3}\right)}{\beta_1 \beta_2 \beta_3^2} d \beta_3 d \beta_2 d \beta_1 d \alpha < 0.00009, \\
\nonumber G^{\prime}_{5} & = 4 \int_{\frac{7}{6}}^{\frac{5}{4}} 1 d \alpha = \frac{1}{3}, \\
\nonumber G^{\prime}_{6} & = \int_{\frac{7}{6}}^{\frac{5}{4}} \int_{\alpha-1}^{\sigma(\alpha)}  \frac{\omega\left(\frac{\alpha-\beta}{\beta}\right)}{\beta^2} d \beta d \alpha > 0.05016, \\
\nonumber G^{\prime}_{7} & = 4 \int_{\frac{5}{4}}^{1.317} 1 d \alpha = 0.268.
\end{align}
By a simple calculation, the value of the right hand side of (7) is less than $0.838x$, and the proof of Theorem~\ref{t2} is completed.

\bibliographystyle{plain}
\bibliography{bib}

\begin{thebibliography}{10}

\bibitem{BD2019}
R.~de~la Bretèche and S.~Drappeau.
\newblock Niveau de répartition des polynômes quadratiques et crible majorant pour les entiers friables.
\newblock {\em J. Eur. Math. Soc.}, 22:1577--1624, 2020.

\bibitem{DI82}
J.-M. Deshouillers and H.~Iwaniec.
\newblock On the greatest prime factor of $n^2 +1$.
\newblock {\em Ann. Inst. Fourier (Grenoble)}, 32(4):1--11, 1982.

\bibitem{EverestHarman}
G.~R. Everest and G.~Harman.
\newblock On primitive divisors of $n^2 +b$.
\newblock In {\em Number Theory and Polynomials}, volume 352 of {\em London Math. Soc. Lecture Note Ser.}, pages 142--154. Cambridge Univ. Press, 2008.

\bibitem{ESTW}
G.~R. Everest, S.~Stevens, D.~Tamsett, and T.~Ward.
\newblock Primes generated by recurrence sequences.
\newblock {\em Amer. Math. Monthly}, 114:417--431, 2007.

\bibitem{GrimmeltMerikoski}
L.~Grimmelt and J.~Merikoski.
\newblock {On the greatest prime factor and uniform equidistribution of quadratic polynomials}.
\newblock {\em arXiv e-prints}, page arXiv:2505.00493v1, 2025.

\bibitem{Harman2024}
G.~Harman.
\newblock Two problems on the greatest prime factor of $n^2 +1$.
\newblock {\em Acta Arith.}, 213(3):273--287, 2024.

\bibitem{Hooley}
C.~Hooley.
\newblock On the greatest prime factor of a quadratic polynomial.
\newblock {\em Acta Math.}, 117:281--299, 1967.

\bibitem{IwaniecP2}
H.~Iwaniec.
\newblock Almost--primes represented by quadratic polynomials.
\newblock {\em Invent. Math.}, 47:171--188, 1978.

\bibitem{KimSarnak}
H.~H. Kim.
\newblock Functoriality for the exterior square of {$G L_4$} and the symmetric fourth of {$G L_2$}, with appendix 1 by {D}. {R}amakrishnan and appendix 2 by {H}. {H}. {K}im and {P}. {S}arnak.
\newblock {\em J. Amer. Math. Soc.}, 16:139--183, 2003.

\bibitem{LRBPrimitive}
R.~Li.
\newblock {On the primitive divisors of quadratic polynomials}.
\newblock {\em arXiv e-prints}, page arXiv:2406.07575v1, 2024.

\bibitem{Merikoski}
J.~Merikoski.
\newblock On the greatest prime factor of $n^2 +1$.
\newblock {\em J. Eur. Math. Soc.}, 25:1253--1284, 2023.

\bibitem{Mine}
M.~Mine.
\newblock {An upper bound for the number of smooth values of a polynomial and its applications}.
\newblock {\em arXiv e-prints}, page arXiv:2410.09558v1, 2024.

\bibitem{Pascadi}
A.~{Pascadi}.
\newblock {Large sieve inequalities for exceptional Maass forms and applications}.
\newblock {\em arXiv e-prints}, page arXiv:2404.04239v1, 2024.

\bibitem{RichertP3}
H.-E. Richert.
\newblock Selberg's sieve with weights.
\newblock {\em Mathmatika}, 16:1--22, 1969.

\bibitem{Schinzel}
A.~Schinzel.
\newblock On two theorems of {G}elfond and some of their applications.
\newblock {\em Acta Arith.}, 13:177--236, 1967.

\bibitem{Tenenbaum}
G.~Tenenbaum.
\newblock {\em Introduction to analytic and probabilistic number theory}, volume 163 of {\em Graduate Studies in Mathematics}.
\newblock American Mathematical Society, Providence, RI, third edition, 2015.

\end{thebibliography}
\end{document}